\documentclass[letterpaper,12pt]{paper}
\usepackage{amsmath,amssymb,amsthm,fullpage}
\usepackage{mathtext} % если нужны русские буквы в формулах (не обязате
\usepackage[cp1251]{inputenc} % кодировка - [koi8-r][cp866][utf8][cp1251]
\usepackage[T2A]{fontenc}  % внутренняя T2A кодировка TeX
\usepackage[english]{babel}   % включение переносов
\usepackage{cite}
\usepackage[dvips]{graphicx}
\usepackage{mathrsfs}
\usepackage{hyperref}

\numberwithin{equation}{section}

\renewcommand\le{\leqslant}
\renewcommand\ge{\geqslant}

\newcommand{\Rset}{\mathbb{R}}
\newcommand{\Cset}{\mathbb{C}}

\renewcommand{\Re}{\mathop{\mathrm{Re}}\nolimits}
\renewcommand{\Im}{\mathop{\mathrm{Im}}\nolimits}

\def\Ai{\mathrm{Ai}}
\def\Bi{\mathrm{Bi}}

\newcommand{\Beta}{\mathrm{B}}
\newcommand{\Mgoth}{\mathfrak{M}}

\newcommand{\Fscr}{\mathscr{F}}
\def\Rint2{\int\!\!\!\int_{\Rset^2}}

\renewcommand\le{\leqslant}
\renewcommand\ge{\geqslant}

\newcommand{\lshad}{[\![}
\newcommand{\rshad}{]\!]}
\newcommand{\rme}{\mathrm{e}}
\newcommand{\rmi}{\mathrm{i}}
\newcommand{\rmd}{\mathrm{d}}

\newcommand{\sgn}{\mathop{\mathrm{sgn}}\nolimits}

\title{Elliptic umbilic representations\\ connected with the caustic}

\author{E.G. Abramochkin and E.V. Razueva\\ 
Coherent Optics Lab, Lebedev Physical Institute,\\ 
Samara, 443011, Russia\\
}
\begin{document}

\maketitle
\begin{abstract}
We investigate the elliptic umbilic canonical integral with an approach based 
on a series expansion of its initial distribution shifted to the caustic 
points. An absolutely convergent integral representation for the elliptic 
umbilic is obtained. Using it, we find the elliptic umbilic particular values 
in terms of 2F2 hypergeometric functions. We also derive an integral over the 
product of Gaussian and two Airy functions in terms of Bessel functions of 
fractional orders. Some other corollaries including 3F2 hypergeometric 
function special values and the Airy polynomials relations are also discussed.
\end{abstract}


\section{Introduction}\label{sec:Intro}
%%%%%%%%%%%%%%%%%%%%%%%%%%%%%%%%%%%%%%%%%%%%%%%%%%%%%%%%%%%%%%%%%%%%%%%%%%%%

This paper is inspired by the relation connecting the elliptic umbilic 
integral \cite{Berry1979} with Bessel functions of orders $\pm 1/6$, which 
was found by Berry and Howls in \cite{Berry2010}:
\begin{align}
E(0,0,z)&=\sqrt{\frac{\pi|z|}{27}}\,\exp\biggl(-\frac{2\rmi z^3}{27}\biggr)
\biggl\{J_{-1/6}\biggl(\frac{2|z|^3}{27}\biggr)
-\rmi\sgn z\cdot J_{1/6}\biggl(\frac{2|z|^3}{27}\biggr)\!\biggr\}
\label{eq:Umbilic_00z}\\
&{}=\frac{\sqrt\pi}{3}\exp\biggl(-\frac{2\rmi z^3}{27}\biggr)
\biggl\{\frac{1}{\Gamma\bigl(\tfrac56\bigr)}\cdot
{}_0F_1\biggl(\frac56\,\Bigl|\,-\frac{z^6}{3^6}\biggr)
-\frac{\rmi z}{3\Gamma\bigl(\tfrac76\bigr)}\cdot
{}_0F_1\biggl(\frac76\,\Bigl|\,-\frac{z^6}{3^6}\biggr)\!\biggr\}.
\label{eq:Umbilic_HyperG_00z}
\end{align}
(We prefer the hypergeometric description of the result since it demonstrates 
that $E(0,0,z)$ is an entire function of $z$. Besides, it leads to the value 
$E(0,0,0)$ immediately. On the other hand, the expression in terms of Bessel 
functions is useful for large values of $z$ due to well-known asymptotic of 
$J_\nu(z)$.) The elliptic umbilic integral is defined by the formula
\begin{equation}
E(x,y,z)=\Fscr\bigl[\exp\bigl\{-\rmi z(\xi^2+\eta^2)
+\rmi(\xi^3-3\xi\eta^2)\bigr\}\bigr](x,y),
\label{eq:Umbilic}
\end{equation}
where $\Fscr[f(\xi,\eta)](x,y)$ is the 2D Fourier transform,
\begin{equation}
\Fscr[f(\xi,\eta)](x,y)=\frac{1}{2\pi}
\Rint2 \exp\bigl\{-\rmi(x\xi+y\eta)\bigr\}f(\xi,\eta)\,\rmd\xi\,\rmd\eta,
\label{eq:Fourier}
\end{equation}
and $x$, $y$, $z$ are real variables. (We use the notation proposed in 
\cite{Berry1979} instead of its later modification 
$\Psi_{\mathrm E}(x,y,z)=2\pi E(-x,-y,-z)$, see \cite{Berry2010} and 
\cite[Chap.\,36]{DLMF}.)

It was quite unexpected for us that $E(0,0,z)$ requires the use of 
$J_{\pm 1/6}$ functions, especially since $E(x,y,0)$ is described in terms of 
Airy functions, i.e., in terms of $J_{\pm 1/3}$ functions 
\cite[Eq.\,(3.4)]{Berry1979}:
\begin{equation}
E(x,y,0)=\gamma\pi\Re\biggl\{\Ai\biggl(-\frac{x+\rmi y}{3\gamma}\biggr)
\Bi\biggl(-\frac{x-\rmi y}{3\gamma}\biggr)\!\biggr\},
\label{eq:Umbilic_xy0}
\end{equation}
where $\gamma=(2/3)^{2/3}$ is an auxiliary constant used for brevity.

The function $E(x,y,z)$ has been thoroughly investigated by Berry 
\textit{et al} \cite{Berry1979, Berry2010}. We recall some properties of 
$E(x,y,z)$ proven in \cite{Berry1979} and used below. First, the function 
$E(x,y,z)$ satisfies the symmetry relations:
\begin{align}
E(x,y,-z)&=E^*(x,y,z),
\label{eq:Umbilic_Conjugate}\\
E(x,-y,z)&=E(x,y,z),
\label{eq:Umbilic_Even_y}\\
E\biggl(x\cos\frac{2\pi n}{3}-y\sin\frac{2\pi n}{3},\,& 
x\sin\frac{2\pi n}{3}+y\cos\frac{2\pi n}{3},z\biggr)=E(x,y,z),
\label{eq:Umbilic_Rotation}
\end{align}
where an asterisk means complex conjugation, $n$ is an integer, and the last 
expression demonstrates the invariance of $E(x,y,z)$ under rotation by 
$2\pi n/3$ in the $(x,y)$ plane.

Second, the equation of the caustic in parametric form is
\begin{equation}
x+\rmi y=\frac{z^2}{3}\rme^{\rmi t}(2+\rme^{-3\rmi t}),\qquad t\in[0,2\pi).
\label{eq:Umbilic_Caustic}
\end{equation}

And third, $E(x,y,z)$ may be expanded into a series in powers of $z$:
\begin{equation}
E(x,y,z)=\gamma\pi\,\sum_{n\ge 0} \frac{(\rmi\gamma z)^n}{n!}
\Re\biggl\{\Ai^{(n)}\biggl(-\frac{x+\rmi y}{3\gamma}\biggr)
\Bi^{(n)}\biggl(-\frac{x-\rmi y}{3\gamma}\biggr)\!\biggr\},
\label{eq:Umbilic_Series0}
\end{equation}
which provides a simple description of the function behaviour for small $z$ 
(see figure~1).
% for small values of~$z$ (see figure~1). --- чтобы уместилось в одну строку

\begin{figure}[t]
\centerline{\includegraphics[width=0.8\textwidth,keepaspectratio]{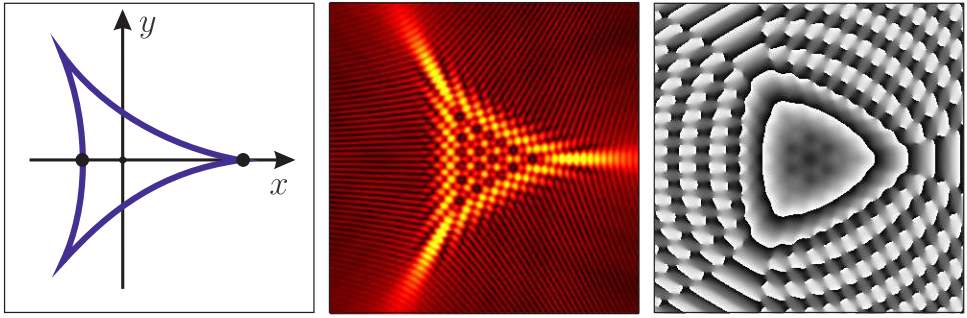}}
\caption[]{\small The caustic, amplitude and phase distributions of the 
elliptic umbilic $E(x,y,z)$. The points $(z^2,0)$ and $(-z^2/3,0)$ marked by 
black discs correspond to the parameter values $t=0$ and $t=\pi$, 
respectively. The distributions are shown for $z=4.5$ in the square 
$-25\le x,y\le 25$.\par}
\label{fig:Fig1}
\end{figure}

Since both functions $\Ai(x)$ and $\Bi(x)$ satisfy the same differential 
equation $y''=xy$, the $n$-th derivative of each function can be written as 
a linear combination of the function itself and its first derivative:
\begin{equation}
\begin{aligned}
\Ai^{(n)}(x)&=P_n(x)\Ai(x)+Q_n(x)\Ai'(x),\\
\Bi^{(n)}(x)&=P_n(x)\Bi(x)+Q_n(x)\Bi'(x).
\end{aligned}
\label{eq:AiryDeriv}
\end{equation}
Here $P_n(x)$ and $Q_n(x)$ are some polynomials (the index $n$ corresponds to 
the derivative order but not the polynomials degree). In particular,
$$
\begin{aligned}
P_0(x)&=1, & P_1(x)&=0, & P_2(x)&=x, & P_6(x)&=x^3+4, & P_{12}(x)&=x^6+260x^3+280,\\
Q_0(x)&=0, & Q_1(x)&=1, & Q_2(x)&=0, & Q_6(x)&=6x,    & Q_{12}(x)&=30x^4+600x.
\end{aligned}
$$
The polynomials can be expressed in terms of particular values of Gegenbauer 
polynomials \cite{Sigma}, however, corresponding formulae are quite unusual 
to apply them for evaluation of the series (\ref{eq:Umbilic_Series0}) easily.

Nevertheless, an attempt to transform and reorder the terms in the series 
expansion (\ref{eq:Umbilic_Series0}) is one of possible ways to find an 
expression of $E(x,y,z)$, which can be helpful for numerical simulations.
Of course, there are various methods for evaluating infinite range oscillatory 
integrals. In particular, the Weniger transform demonstrates the effectiveness
being applied to the elliptic umbilic \cite{Borghi2010, Borghi2016}. However, 
this method is mainly due to divergent integrals and series, whereas the right 
side of (\ref{eq:Umbilic}) is an ordinary integral.
In a sense, the situation with the Airy function, $\Ai(x)$, is a good example. 
Initially, $\Ai(x)$ is defined by a conditionally convergent integral for 
real values of~$x$ only. Then, by contour integration, it reduces to an 
absolutely convergent integral which is valid for any complex~$x$ (see 
Eqs.\,(9.5.1) and (9.5.4) in \cite[Chap.\,9]{DLMF}). Here we are interested in 
similar procedure for the elliptic umbilic.

The paper is organized as follows. In Section~\ref{sec:IntegAiry} we find 
various series expansions for the elliptic umbilic depending on the 
derivarives of the Airy functions product. This section also gives the main 
result of the paper, namely an absolutely convergent integral representation 
of $E(x,y,z)$ closely connecting with the caustic (\ref{eq:Umbilic_Caustic}). 
In Sections~\ref{sec:HyperG_Series} and \ref{sec:Integ_AiBiGauss} some 
corollaries of this result are presented. First, we analyze the series 
expansions and obtain the elliptic umbilic values in the caustic points 
$(z^2,0,z)$ and $(-z^2/3,0,z)$ in terms of ${}_2F_2$ hypergeometric functions. 
Then, an improper integral of the product of Gaussian and two Airy functions 
is evaluated in terms of ${}_0F_1$ hypergeometric functions (the Bessel 
functions of fractional orders). The paper ends with concluding remarks and 
perspectives for future research. Some auxiliary integrals and series 
containing the Airy functions are presented in Appendix~\ref{sec:AppA}.

\section{Elliptic umbilic and absolutely convergent integral representation}
\label{sec:IntegAiry}
%%%%%%%%%%%%%%%%%%%%%%%%%%%%%%%%%%%%%%%%%%%%%%%%%%%%%%%%%%%%%%%%%%%%%%%%%%%%

Let us return to the elliptic umbilic definition (\ref{eq:Umbilic}). 
By changing variables $\xi\to\xi+a$ and $\eta\to\eta+b$, where $a,b$ are real 
parameters, we obtain
\begin{align}
E(x,y,z)&=\frac{1}{2\pi}
\exp\bigl(-\rmi ax-\rmi by-\rmi z[a^2+b^2]+\rmi[a^3-3ab^2]\bigr)
\nonumber\\
&\hspace{10mm}{}
\times\Rint2 \exp\bigl(-\rmi[x+2az-3(a^2-b^2)]\xi-\rmi[y+2bz+6ab]\eta
\nonumber\\
&\hspace{20mm}{}-\rmi[z-3a]\xi^2-6\rmi b\xi\eta-\rmi[z+3a]\eta^2
+\rmi[\xi^3-3\xi\eta^2]\bigr)\,\rmd\xi\,\rmd\eta.
\label{eq:Umbilic_xyz_ab}
\end{align}
Expanding the factor $\exp(-\rmi[z-3a]\xi^2-6\rmi b\xi\eta-\rmi[z+3a]\eta^2)$ 
as a power series and replacing monomials in $\xi,\eta$ by derivatives in 
$x,y$ leads to an integral which is known due to (\ref{eq:Umbilic_xy0}):
%
% Если будет бухтёж рецензента про расходимость интеграла с мономом, то 
% дать вывод через контур, который использовал Берри в \cite{Berry1979} и 
% процитировать Литлвуда (стр.90): ... строгость в анализе за пределами 
% студенческих работ не имеет первостепенного значения и может быть при 
% наличии настоящей идеи всегда внесена любым компетентным профессионалом. 
% ... rigour, which for that matter is not of first-rate importance in 
% analysis beyond the undergraduate stage, and can be supplied, 
% given a real idea, by any competent professional.
% J.E.Littlewood. "A Mathematician's Miscellany" London, Methuen, 1953.
%
\begin{align}
&E(x,y,z)=\exp\bigl(-\rmi ax-\rmi by-\rmi z[a^2+b^2]+\rmi[a^3-3ab^2]\bigr)
\nonumber\\
&{}\times\sum_{n\ge 0} \,\frac{\rmi^n}{n!}
\bigl\{[z-3a]\partial_x^2+6b\partial_x\partial_y+[z+3a]\partial_y^2\bigr\}^n 
E(x+2az-3[a^2-b^2],y+2bz+6ab,0)
\label{eq:Umbilic_Series_ab}
\end{align}
(cf. \cite[Eq.\,(27)]{Berry1980}).

Choosing various values of $a,b$ in (\ref{eq:Umbilic_Series_ab}) gives various 
series expansions. For example, if the elliptic umbilic in the right side of 
(\ref{eq:Umbilic_Series_ab}) is $E(x,y,0)$, then we have 
$$
\begin{array}{c} 2az-3[a^2-b^2]=0\\ 2bz+6ab=0\end{array}\biggr\}
\quad\Rightarrow\quad
(a,b)=\biggl\{(0,0),\,
\biggl(-\frac{z}{3},\frac{z}{\sqrt 3}\biggr),\, 
\biggl(-\frac{z}{3},-\frac{z}{\sqrt 3}\biggr),\,
\biggl(\frac{2z}{3},0\biggr)\!\biggr\}.
% \label{eq:System_ab}
$$
For the case $a=b=0$, the expansion (\ref{eq:Umbilic_Series_ab}) reduces 
to (\ref{eq:Umbilic_Series0}):
\begin{align}
E(x,y,z)&=\sum_{n\ge 0} \,\frac{(\rmi z)^n}{n!}\,
\{\partial_x^2+\partial_y^2\}^n E(x,y,0)
\nonumber\\
&{}=\sum_{n\ge 0} \,\frac{(\rmi z)^n}{n!}\,
\{\gamma\partial_v\partial_{v^*}\}^n E(x,y,0)
\nonumber\\
&{}=\gamma\pi\sum_{n\ge 0} \,\frac{(\rmi\gamma z)^n}{n!}\,
\Re\{\Ai^{(n)}(v)\Bi^{(n)}(v^*)\},
\label{eq:Umbilic_Series00}
\end{align}
where $v=-(x+\rmi y)/3\gamma$ \cite{GKP_note}. For other three cases, the 
expansions are
\begin{align}
E(x,y,z)
&{}=\exp\biggl(\frac{\rmi z}{3}(x-y\sqrt{3})-\frac{4\rmi z^3}{27}\biggr)
\sum_{n\ge 0} \,\frac{(2\rmi z)^n}{n!}\,
\bigl\{\partial_x^2+\sqrt{3}\,\partial_x\partial_y\bigr\}^n E(x,y,0),
\nonumber\\
&{}=\exp\biggl(\frac{\rmi z}{3}(x+y\sqrt{3})-\frac{4\rmi z^3}{27}\biggr)
\sum_{n\ge 0} \,\frac{(2\rmi z)^n}{n!}\,
\bigl\{\partial_x^2-\sqrt{3}\,\partial_x\partial_y\bigr\}^n E(x,y,0),
\nonumber\\
&{}=\exp\biggl(-\frac{2\rmi xz}{3}-\frac{4\rmi z^3}{27}\biggr)
\sum_{n\ge 0} \,\frac{(\rmi z)^n}{n!}\,
\bigl\{3\partial_y^2-\partial_x^2\bigr\}^n E(x,y,0).
\label{eq:Umbilic5}
\end{align}
It is evident that all of them are various versions of one and the same 
expansion due to the symmetry relation (\ref{eq:Umbilic_Rotation}).

We can simplify the differential operator in the right side of 
(\ref{eq:Umbilic_Series_ab}) by removing two terms of three:
\begin{align}
E(x,y,z)
&{}=\exp\biggl(\frac{\rmi xz}{3}-\frac{4\rmi z^3}{27}\biggr)
\sum_{n\ge 0} \,\frac{(2\rmi z)^n}{n!}\,\partial_x^{2n}E(x-z^2,y,0),
\label{eq:Umbilic_Series-z/3}\\
&{}=\exp\biggl(-\frac{\rmi xz}{3}-\frac{2\rmi z^3}{27}\biggr)
\sum_{n\ge 0} \,\frac{(2\rmi z)^n}{n!}\,
\partial_y^{2n}E\biggl(x+\frac{z^2}{3},y,0\biggr).
\label{eq:Umbilic_Series+z/3}
\end{align}
Here $(a,b)=(-z/3,0)$ and $(a,b)=(z/3,0)$, respectively. Moreover, expanding 
$E(x,y,0)$ as a power series in $x,y$ (see Appendix~A), it is easy to find the 
derivatives in both expressions.

And finally, substituting $z=0$ into (\ref{eq:Umbilic_Series_ab}) one can 
find a differential difference relation for $E(x,y,0)$:
\begin{align}
E(x,y,0)&=\exp\bigl(-\rmi[ax+by]+\rmi[a^3-3ab^2]\bigr)
\nonumber\\
&{}\times\sum_{n\ge 0} \,\frac{(3\rmi)^n}{n!}
\bigl\{a[\partial_y^2-\partial_x^2]+2b\partial_x\partial_y\bigr\}^n 
E(x-3[a^2-b^2],y+6ab,0).
\label{eq:Umbilic4}
\end{align}

The formulae above are useful in certain analytic descriptions (special values 
of the hypergeometric ${}_3F_2$ function, Airy polynomials $\mathrm{Pi}_n$, 
etc.) which we discuss later. However, all the expansions of $E(x,y,z)$ 
did not look perfect because their invariance under a rotation by $2\pi/3$ 
in $(x,y)$ plane cannot be read immediately. The expansion 
(\ref{eq:Umbilic_Series00}) looks a little better than others since its 
invariance follows from the known properties of the Airy functions:
\begin{gather}
\Ai(\omega v)=\tfrac12\rme^{\pi\rmi/3}\bigl[\Ai(v)-\rmi\cdot\Bi(v)\bigr],\quad
\Bi(\omega v)=\tfrac12\rme^{-\pi\rmi/6}\bigl[3\Ai(v)+\rmi\cdot\Bi(v)\bigr],
\label{eq:AlgRels}
\end{gather}
where $\omega=\rme^{2\pi\rmi/3}$ is the cubic root of unity.

We can try to construct the desired expansion, noting that the elliptic 
umbilic shifts in the formulae (\ref{eq:Umbilic_Series-z/3}) and 
(\ref{eq:Umbilic_Series+z/3}) are coordinates of the caustic points placed on 
the $x$ axis (see figure~1). Namely, we select parameters $a,b$ in 
(\ref{eq:Umbilic_Series_ab}) such that the elliptic umbilic in the right side 
has been shifted to a point of the caustic, $E(x-x_0,y-y_0,0)$, where
\begin{equation}
x_0=\frac{z^2}{3}(2\cos t_0+\cos 2t_0),\quad
y_0=\frac{z^2}{3}(2\sin t_0-\sin 2t_0)
\label{eq:CausticPoint}
\end{equation}
for some $t_0\in[0,2\pi)$. Then
\begin{gather}
\begin{array}{c} 2az-3[a^2-b^2]=-x_0\\ 2bz+6ab=-y_0\end{array}\biggr\}
\quad\Rightarrow\quad
a=-\frac{z}{3}\cos t_0,\quad b=-\frac{z}{3}\sin t_0,
\label{eq:CausticValues_ab}\\
\bigl\{[z-3a]\partial_x^2+6b\partial_x\partial_y+[z+3a]\partial_y^2\bigr\}^n 
=(2z)^n\Bigl\{\cos\frac{t_0}{2}\,\partial_x
-\sin\frac{t_0}{2}\,\partial_y\Bigr\}^{2n},
\nonumber 
% \label{eq:DiffOperator_2n}
\end{gather}
and the expansion (\ref{eq:Umbilic_Series_ab}) has the form
\begin{align}
E(x,y,z)&=\exp\biggl(\frac{\rmi z}{3}[x\cos t_0+y\sin t_0]
-\frac{\rmi z^3}{27}[3+\cos 3t_0]\biggr)
\nonumber\\
&{}\times\sum_{n\ge 0} \frac{(2\rmi z)^n}{n!}
\biggl\{\cos\frac{t_0}{2}\,\partial_x
-\sin\frac{t_0}{2}\,\partial_y\biggr\}^{2n}E(x-x_0,y-y_0,0).
\label{eq:Umbilic_ab2}
\end{align}

Now, we transform the right side of (\ref{eq:Umbilic_ab2}) to an integral 
with Airy functions. Let $X=x-x_0$, $Y=y-y_0$ and $V\!=-(X+\rmi Y)/3\gamma$. 
Then $E(X,Y,0)=\gamma\pi\Re\bigl\{\Ai(V)\Bi(V^*)\bigr\}$,
$$
\cos\frac{t_0}{2}\,\partial_X-\sin\frac{t_0}{2}\,\partial_Y=-\frac{1}{3\gamma}
\bigl\{\rme^{-\rmi t_0/2}\partial_V+\rme^{\rmi t_0/2}\partial_{V^*}\bigr\}
% \label{eq:DiffOperator_1}
$$
and
\begin{align}
&\biggl\{\cos\frac{t_0}{2}\,\partial_X
-\sin\frac{t_0}{2}\,\partial_Y\biggr\}^{\!2n}E(X,Y,0)={}
\nonumber\\
&{}=\frac{\gamma\pi}{2}\biggl(-\frac{1}{3\gamma}\biggr)^{\!2n}
\sum_{k=0}^{2n} \binom{2n}{k}
\bigl(\rme^{-\rmi t_0/2}\partial_V)^k
\bigl(\rme^{\rmi t_0/2}\partial_{V^*}\bigr)^{2n-k}
\bigl\{\Ai(V)\Bi(V^*)+\Ai(V^*)\Bi(V)\bigr\}={}
\nonumber\\
&{}=\frac{\gamma\pi}{2}\biggl(\frac{\gamma}{4}\biggr)^{\!n}
\sum_{k=0}^{2n} \binom{2n}{k}
\bigl[\rme^{\rmi t_0(n-k)}\Ai^{(k)}(V)\Bi^{(2n-k)}(V^*)
+\rme^{-\rmi t_0(n-k)}\Ai^{(k)}(V^*)\Bi^{(2n-k)}(V)\bigr].
\label{eq:DiffOperator_2}
\end{align}
Denoting by $S(n,k)$ the expression within the square brackets and 
substituting the right side of (\ref{eq:DiffOperator_2}) into the series 
in (\ref{eq:Umbilic_ab2}), we have
\begin{align}
{\textstyle \sum\limits_n}&\overset{\textrm{def}}{=}
\sum_{n\ge 0} \frac{(2\rmi z)^n}{n!}\biggl\{\cos\frac{t_0}{2}\,\partial_X
-\sin\frac{t_0}{2}\,\partial_Y\biggr\}^{\!2n}E(X,Y,0)={}
\nonumber\\
&{}=\frac{\gamma\pi}{2}\sum_{n\ge 0} \frac{1}{n!}
\biggl(\frac{\rmi\gamma z}{2}\biggr)^{\!n}\sum_{k=0}^{2n} \binom{2n}{k}S(n,k)
=\frac{\gamma\pi}{2}\sum_{k\ge 0}\sum_{n\ge\lceil\frac{k}{2}\rceil} 
\!\!\binom{2n}{k}\frac{1}{n!}\biggl(\frac{\rmi\gamma z}{2}\biggr)^{\!n}
S(n,k)={}
\nonumber\\
&{}=\frac{\gamma\pi}{2}\sum_{\epsilon=0,1}\sum_{k\ge 0}
\sum_{n\ge k+\epsilon} \!\binom{2n}{2k+\epsilon}
\frac{1}{n!}\biggl(\frac{\rmi\gamma z}{2}\biggr)^{\!n}S(n,2k+\epsilon)={}
\nonumber\\
&{}=\frac{\gamma\pi}{2}\sum_{\epsilon=0,1}\sum_{k\ge 0}
\sum_{n\ge 0} \frac{(2n+2k+2\epsilon)!}{(n+k+\epsilon)!}
\biggl(\frac{\rmi\gamma z}{2}\biggr)^{\!n+k+\epsilon}
\frac{S(n+k+\epsilon,2k+\epsilon)}{(2k+\epsilon)!\,(2n+\epsilon)!}.
\label{eq:Series1}
\end{align}
Since
$$
\frac{(2N)!}{N!}=\frac{1}{\sqrt\pi}\int_\Rset \rme^{-t^2}(4t^2)^N\,\rmd t,
$$
then
\begin{align}
{\textstyle \sum\limits_n}&=\frac{\gamma\sqrt\pi}{2}\sum_{\epsilon=0,1}
\int_\Rset \rme^{-t^2} \sum_{k\ge 0}\sum_{n\ge 0} 
\frac{(2\rmi\gamma zt^2)^{n+k+\epsilon}}{(2k+\epsilon)!\,(2n+\epsilon)!}
\nonumber\\
&{}\times
\bigl\{\rme^{\rmi t_0(n-k)}\Ai^{(2k+\epsilon)}(V)\Bi^{(2n+\epsilon)}(V^*)
+\rme^{-\rmi t_0(n-k)}\Ai^{(2k+\epsilon)}(V^*)\Bi^{(2n+\epsilon)}(V)\bigr\}
\,\rmd t.
\label{eq:Series2}
\end{align}
Returning to (\ref{eq:Umbilic_ab2}) and using the relations
\begin{align}
\sum_{k\ge 0} F^{(2k)}(x)\frac{a^{2k}}{(2k)!}
=\sum_{k\ge 0} \frac{1+(-1)^k}{2}\cdot F^{(k)}(x)\frac{a^k}{k!}
&=\frac12\bigl\{F(x+a)+F(x-a)\bigr\},
\nonumber\\
\sum_{k\ge 0} F^{(2k+1)}(x)\frac{a^{2k+1}}{(2k+1)!}
=\sum_{k\ge 0} \frac{1-(-1)^k}{2}\cdot F^{(k)}(x)\frac{a^k}{k!}
&=\frac12\bigl\{F(x+a)-F(x-a)\bigr\},
\nonumber
\end{align}
we obtain an absolutely convergent integral because both Airy functions are of 
order $\tfrac32$:
\begin{align}
E(x,y,z)&=\frac{\gamma\sqrt{\pi}}{2}
\exp\biggl(\frac{\rmi z}{3}[x\cos t_0+y\sin t_0]
-\frac{\rmi z^3}{27}[3+\cos 3t_0]\biggr)
\nonumber\\
&\times\int_\Rset \rme^{-t^2}
\Bigl\{\Ai(V+t\sqrt{2\rmi\rme^{-\rmi t_0}\gamma z}\,)
\Bi(V^*+t\sqrt{2\rmi\rme^{\rmi t_0}\gamma z}\,)+{}
\nonumber\\
&\hspace{20mm}{}+\Ai(V^*+t\sqrt{2\rmi\rme^{\rmi t_0}\gamma z}\,)
\Bi(V+t\sqrt{2\rmi\rme^{-\rmi t_0}\gamma z}\,)\!\Bigr\}\,\rmd t.
\label{eq:Umbilic_ab4}
\end{align}
Here, $V=-\{(x-x_0)+\rmi(y-y_0)\}/3\gamma$ and $(x_0,y_0)$ is a point of the 
caustic (\ref{eq:CausticPoint}). In particular, if $z=0$, then $x_0=y_0=0$ 
and (\ref{eq:Umbilic_xy0}) follows immediately. As is seen from 
(\ref{eq:Series2}), the expression (\ref{eq:Umbilic_ab4}) does not depend on 
a square root branch selection.

\section{Elliptic umbilic and hypergeometric functions}
\label{sec:HyperG_Series}
%%%%%%%%%%%%%%%%%%%%%%%%%%%%%%%%%%%%%%%%%%%%%%%%%%%%%%%%%%%%%%%%%%%%%%%%%%%%

In what follows, we will restrict ourselves to the case $z\ge 0$ for 
simplicity (see Eq.~(\ref{eq:Umbilic_Conjugate})), considering $z$ as a 
parameter. To evaluate the integral (\ref{eq:Umbilic_ab4}), a caustic point 
$(x_0,y_0)$ is required. Of course, this point may be chosen arbitrary. 
However, for evaluation of $E(x,y,z)$ it seems natural to choose such 
$(x_0,y_0)$ that $x=cx_0$, $y=cy_0$, where $c\ge 0$. In particular, 
if $(x,y)$ is placed on the caustic, the choice $x_0=x$, $y_0=y$ is the 
best since then $V=0$ and the oscillating influence of the Airy functions 
on the integrand is minimal. Nevertheless, we will not exclude the case $c<0$ 
hoping to obtain a nice formula.

Substituting $x=cx_0$ and $y=cy_0$ into (\ref{eq:Umbilic_ab4}), one gets
\begin{align}
E(cx_0,cy_0,z)&=\frac{\gamma\sqrt{\pi}}{2}
\exp\biggl(\frac{\rmi cz^3}{9}[2+\cos 3t_0]
-\frac{\rmi z^3}{27}[3+\cos 3t_0]\biggr)
\nonumber\\
&{}\times\!\int_\Rset \rme^{-t^2}
\bigl\{\Ai(X+\rmi Y)\Bi(X-\rmi Y)+\Ai(X-\rmi Y)\Bi(X+\rmi Y)\bigr\}\,\rmd t,
\label{eq:Umbilic_ab5}
\end{align}
where
\begin{equation}
\begin{aligned}
X&=\frac{V+V^*}{2}+t\sqrt{2\rmi\gamma z}\,
\frac{\rme^{-\rmi t_0/2}+\rme^{\rmi t_0/2}}{2}
=\frac{(1-c)z^2}{9\gamma}(2\cos t_0+\cos 2t_0)
+t\sqrt{2\rmi\gamma z}\,\cos\frac{t_0}{2},\\
Y&=\frac{V-V^*}{2\rmi}+t\sqrt{2\rmi\gamma z}\,
\frac{\rme^{-\rmi t_0/2}-\rme^{\rmi t_0/2}}{2\rmi}
=\frac{(1-c)z^2}{9\gamma}(2\sin t_0-\sin 2t_0)
-t\sqrt{2\rmi\gamma z}\,\sin\frac{t_0}{2}.
\end{aligned}
\label{eq:Values_XY}
\end{equation}
Applying the results of Appendix~A, we can expand the expression within the 
figure brackets into a Taylor series:
\begin{align}
\Ai(X+\rmi Y)\Bi(X-\rmi Y)&{}+\Ai(X-\rmi Y)\Bi(X+\rmi Y)={}
\nonumber\\
&{}=\frac{8}{\sqrt\pi}
\sum_{m,n\ge 0} \,\frac{(-1)^m\{\Im\omega^{n-m+1}\} X^m Y^{2n}}
{12^{(5+2n-2m)/6}\Gamma\bigl(\tfrac{5+2n-2m}{6}\bigr)\,m!\,n!}.
\label{eq:AiBi+AiBi}
\end{align}
Then substituting (\ref{eq:Values_XY}) into (\ref{eq:AiBi+AiBi}), expanding 
by Newton's binomial formula and integrating term by term converts the 
integral in (\ref{eq:Umbilic_ab5}) into a multiple series:
\begin{align}
E(cx_0,cy_0,z)&=\frac{2\sqrt\pi}{3\sqrt 3}
\exp\biggl(\bigl\{[6c-3]+[3c-1]\cos 3t_0\bigr\}\frac{\rmi z^3}{27}\biggr)
\sum_{m,n\ge 0} \frac{\{\Im\omega^{n-m+1}\}}
{\Gamma\bigl(\tfrac{5+2n-2m}{6}\bigr)m!\,n!}
\nonumber\\
&{}\times\sum_{0\le 2\ell\le m+2n} T(\ell,m,n)\biggl(\frac12\biggr)_{\!\ell}
\frac{4^\ell(c-1)^{m+2n-2\ell}}{3^{m+3n-2\ell}}(\rmi z)^{2(m+2n)-3\ell},
\label{eq:Umbilic_ab6}
\end{align}
where
\begin{align}
T(\ell,m,n)&=\sum_{k\ge 0} \frac{1}{2^k}\binom{m}{2\ell-k}\!\binom{2n}{k}
\nonumber\\
&{}\times(1-2\cos t_0-2\cos^2 t_0)^{m-2\ell+k}
(1+\cos t_0)^{n+\ell-k}(1-\cos t_0)^{3n-k}.
\label{eq:Series_T}
\end{align}
The series $T(\ell,m,n)$ is naturally terminating, so we do not indicate its 
limits here. (In fact, $\max(0,2\ell-m)\le k\le\min(2\ell,2n)$ and 
$T(\ell,m,n)$ is a terminating ${}_2F_1$ hypergeometric series.)

The right side of (\ref{eq:Umbilic_ab6}), after some algebra, may be written 
in the manner of (\ref{eq:Umbilic_HyperG_00z}),
\begin{align}
E(cx_0,cy_0,z)&=\frac{\sqrt{\pi}}{3}
\exp\biggl(\bigl\{[6c-3]+[3c-1]\cos 3t_0\bigr\}\frac{\rmi z^3}{27}\biggr)
\nonumber\\
&{}\times\biggl\{\frac{1}{\Gamma\bigl(\tfrac56\bigr)}
\sum_{\nu\ge 0} W_{0,3\nu}\biggl(\frac{\rmi z}{3}\biggr)^{3\nu}
-\frac{\rmi z}{3\Gamma\bigl(\tfrac76\bigr)}
\sum_{\nu\ge 0} W_{1,3\nu+1}\biggl(\frac{\rmi z}{3}\biggr)^{3\nu}\biggr\},
\label{eq:Umbilic_ab7}
\end{align}
where
\begin{align}
W_{\delta,N}&=\sum_{0\le 2k\le N} 
4^{N-2k}\biggl(\frac12\biggr)_{\!N-2k}(c-1)^k\times{}
\nonumber\\
&{}\times\sum_{0\le 2n\le 2N-3k} 3^{N-k-n}\,
\frac{\Gamma\bigl(\tfrac{5+2\delta}{6}\bigr)\,T(N-2k,2N-3k-2n,n)}
{\Gamma\bigl(\tfrac{5-4N}{6}+n+k\bigr)(2N-3k-2n)!\,n!}.
\label{eq:Coeffs_W}
\end{align}
In particular, $W_{0,0}=W_{1,1}=1$, $W_{0,3}=(6c-1)+(3c-1)\cos 3t_0$,
$W_{1,4}=\bigl(\tfrac{15}{4}c^2+6c-1\bigr)+(3c^2+3c-1)\cos 3t_0$. 
Unfortunately, we could not rewrite $W_{\delta,N}$ as polynomials in 
$\cos 3t_0$, in general.

Now, we consider two cases of (\ref{eq:Umbilic_ab6}), when its quadruple 
series can be shortened drastically. These are $t_0=0$ and $t_0=\pi$ that 
correspond to caustic points placed on the $x$ axis:
\begin{align}
E(cz^2,0,z)&=\frac{2\sqrt\pi}{3\sqrt 3}
\exp\biggl([9c-4]\frac{\rmi z^3}{27}\biggr)\sum_{m,n\ge 0} 
\{\Im\omega^{n-m+1}\}\frac{([c-1]z^2)^m(2\rmi z)^n}
{\Gamma\bigl(\tfrac{5-2m-4n}{6}\bigr)m!\,n!},
\label{eq:Umbilic_XRc}\\
E\biggl(-\frac{cz^2}{3},0,z\biggr)&=\frac{2\sqrt\pi}{3\sqrt 3}
\exp\biggl([3c-2]\frac{\rmi z^3}{27}\biggr)\sum_{m,n\ge 0} 
\{\Im\omega^{n-m+1}\}\frac{\bigl(\tfrac12\bigr)_n
\bigl(\tfrac13[1-c]z^2\bigr)^m\bigl(\tfrac23\rmi z\bigr)^n}
{\Gamma\bigl(\tfrac{5-2m+2n}{6}\bigr)m!\,n!}.
\label{eq:Umbilic_XLc}
\end{align}
Since $\omega^{n-m+1}=\omega^{1+2m+n}$, it is reasonable to collect terms 
with $2m+n=k$ together. Then we have
\begin{align}
E(cz^2,0,z)&=\frac{2\sqrt\pi}{3\sqrt 3}
\exp\biggl([9c-4]\frac{\rmi z^3}{27}\biggr)
\nonumber\\
&{}\times\sum_{k\ge 0} \,\{\Im\omega^{k+1}\}\cdot
\frac{(2\rmi z)^k}{\Gamma\bigl(\tfrac{5-4k}{6}\bigr)k!}\cdot
{}_2F_1\Bigl(\frac{-k}{2},\frac{1-k}{2};\,\frac{5-4k}{6}\,\Bigl|\,1-c\Bigr),
\label{eq:Umbilic_XRc2}\\
E\biggl(-\frac{cz^2}{3},0,z\biggr)&=\frac{2\sqrt\pi}{3\sqrt 3}
\exp\biggl([3c-2]\frac{\rmi z^3}{27}\biggr)
\nonumber\\
&{}\times\sum_{k\ge 0} \,\{\Im\omega^{k+1}\}
\frac{\bigl(\tfrac12\bigr)_k\bigl(\tfrac23\rmi z\bigr)^k}
{\Gamma\bigl(\tfrac{5+2k}{6}\bigr)k!}\cdot{}_3F_2\biggl(\!\begin{array}{c} 
\frac{-k}{2},\frac{1-k}{2},\frac{1-2k}{6}\smallskip\\ 
\frac{1-2k}{4},\frac{3-2k}{4} \end{array}\!\Bigl|\,\frac34[1-c]\biggr),
\label{eq:Umbilic_XLc2}
\end{align}
where the ${}_2F_1$ functions may be reduced to Gegebauer polynomials. 

Of course, the case $c=1$ is the simplest for both formulae:
\begin{align}
E(z^2,0,z)&=\frac{\sqrt\pi}{3}\exp\biggl(\frac{5\rmi z^3}{27}\biggr)
\nonumber\\
&{}\times\biggl\{\frac{1}{\Gamma\bigl(\tfrac56\bigr)}\cdot
{}_2F_2\biggl(\!\!\begin{array}{c} \frac{1}{12},\frac{7}{12} \smallskip\\ 
\frac13,\frac23\end{array}\!\Bigl|\,{-}\frac{32}{27}\rmi z^3\!\biggr)
\!-\!\frac{\rmi z}{3\Gamma\bigl(\tfrac76\bigr)}\cdot
{}_2F_2\biggl(\!\!\begin{array}{c} \frac{5}{12},\frac{11}{12} \smallskip\\ 
\frac23,\frac43\end{array}\!\Bigl|\,{-}\frac{32}{27}\rmi z^3\!\biggr)
\!\biggr\},
\label{eq:Umbilic_CausticXR}\\
E\biggl(-\frac{z^2}{3},0,z\biggr)
&=\frac{\sqrt\pi}{3}\exp\biggl(\frac{\rmi z^3}{27}\biggr)
\nonumber\\
&{}\times\biggl\{\frac{1}{\Gamma\bigl(\tfrac56\bigr)}\cdot
{}_2F_2\biggl(\!\!\begin{array}{c} \frac16,\frac12 \smallskip\\ 
\frac13,\frac23\end{array}\!\Bigl|\,{-}\frac{8}{27}\rmi z^3\!\biggr)
\!-\!\frac{\rmi z}{3\Gamma\bigl(\tfrac76\bigr)}\cdot
{}_2F_2\biggl(\!\!\begin{array}{c} \frac12,\frac56 \smallskip\\ 
\frac23,\frac43\end{array}\!\Bigl|\,{-}\frac{8}{27}\rmi z^3\!\biggr)
\!\biggr\}.
\label{eq:Umbilic_CausticXL}
\end{align}

If $c=0$, then (\ref{eq:Umbilic_XRc2}) leads to (\ref{eq:Umbilic_HyperG_00z}) 
due to Gauss formula for a ${}_2F_1$ of unit argument 
\cite[Eq.\,(15.4.20)]{DLMF} 
and the relation between Kummer and Bessel functions, ${}_1F_1(a;\,2a\,|\,4z)
=\rme^{2z}\cdot{}_0F_1\bigl(a+\tfrac12\,\bigl|\,z^2\bigr)$.

The asymptotic expansion of the ${}_2F_2$ function as $z\gg 1$ is well known 
\cite[Sect.\,16.11]{DLMF}. As result,
\begin{align}
E(z^2,0,z)
&\approx\frac{1}{2\sqrt{6\pi}}\exp\Bigl(\frac{5\rmi z^3}{27}\Bigr)
\biggl\{\rme^{-\pi\rmi/8}\,\frac{\Gamma\bigl(\tfrac14\bigr)}{(2z)^{1/4}}
-\rme^{\pi\rmi/8}\,\frac{\Gamma\bigl(\tfrac74\bigr)}{(2z)^{7/4}}\biggr\}
+\frac{\exp(-\rmi z^3)}{4\sqrt{2}\,z},
\label{eq:Umbilic_Asympt_XR}\\
E\biggl(-\frac{z^2}{3},0,z\biggr)
&\approx\frac{\sqrt 6}{2\sqrt\pi}\exp\Bigl(\frac{\rmi z^3}{27}\Bigr)
\biggl\{\rme^{-\pi\rmi/4}\,\frac{\Gamma\bigl(\tfrac43\bigr)}{(4z)^{1/2}}
-\rme^{\pi\rmi/4}\,\frac{\Gamma\bigl(\tfrac23\bigr)}{(4z)^{3/2}}\biggr\}
+\frac{1}{2z}\exp\Bigl(-\frac{7\rmi z^3}{27}\Bigr).
\label{eq:Umbilic_Asympt_XL}
\end{align}
(Each of the terms should be multiplied by $\{1+\mathcal{O}(z^{-3})\}$ to 
see the remainder term order.)

Following \cite{Berry2010}, we compare the asymptotic expansions 
(\ref{eq:Umbilic_Asympt_XR}) and (\ref{eq:Umbilic_Asympt_XL}) with the exact 
solutions (\ref{eq:Umbilic_CausticXR}) and (\ref{eq:Umbilic_CausticXL}) in 
figure~2.

\begin{figure}[t]
\centerline{\includegraphics[width=0.8\textwidth,keepaspectratio]{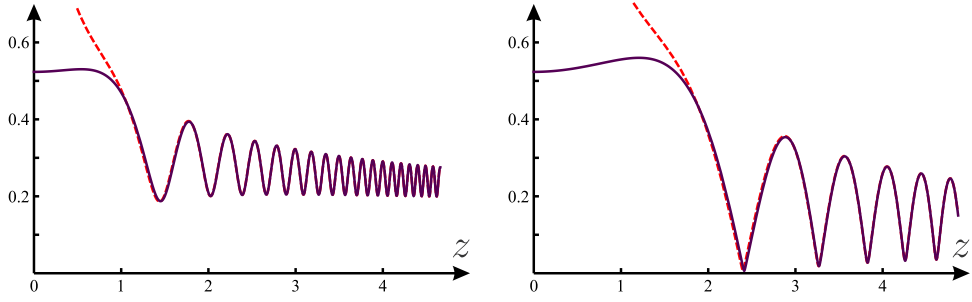}}
\caption[]{\small The absolute value of $E(cz^2,0,z)$ (full curve) and its 
asymptotic approximation (dashed curve) for $c=1$~(left) and 
$c=-1/3$~(right).\par}
\label{fig:Fig2}
\end{figure}

\section{Integrals containing the products of Airy functions}
\label{sec:Integ_AiBiGauss}
%%%%%%%%%%%%%%%%%%%%%%%%%%%%%%%%%%%%%%%%%%%%%%%%%%%%%%%%%%%%%%%%%%%%%%%%%%%%

Let us return to the integral representation (\ref{eq:Umbilic_ab4}) and 
consider the cases, for which the value of $E(x,y,z)$ is already known, 
Eqs.~(\ref{eq:Umbilic_HyperG_00z}), (\ref{eq:Umbilic_CausticXR}) and 
(\ref{eq:Umbilic_CausticXL}). Then the integral in the right side of 
(\ref{eq:Umbilic_ab4}) may be expressed in hypergeometric terms. Moreover, 
since both sides of (\ref{eq:Umbilic_ab4}) are entire functions, it holds for 
the complex values of $z$ and $t_0$ also.

We start with the case $E(0,0,z)$. Introducing parameters 
$2a=\sqrt{2\rmi\rme^{-\rmi t_0}\gamma z}$ and 
$2b=\sqrt{2\rmi\rme^{\rmi t_0}\gamma z}$, we get the relation
\begin{align}
&\int_\Rset \rme^{-t^2}
\Bigl\{\Ai\bigl(2at-a[a^3+2b^3]\bigr)\Bi\bigl(2bt-b[2a^3+b^3]\bigr)
\nonumber\\
&\hspace{50mm}{}+\Ai\bigl(2bt-b[2a^3+b^3]\bigr)
\Bi\bigl(2at-a[a^3+2b^3]\bigr)\!\Bigr\}\,\rmd t
\nonumber\\
&{}=\frac{2}{12^{1/3}}\exp\biggl(-\frac{(a^3+b^3)^2}{3}\biggr)
\biggl\{\frac{1}{\Gamma\bigl(\tfrac56\bigr)}\cdot
{}_0F_1\biggl(\frac56\,\Bigl|\,\frac49(ab)^6\biggr)
-\frac{2ab}{12^{1/3}\Gamma\bigl(\tfrac76\bigr)}\cdot
{}_0F_1\biggl(\frac76\,\Bigl|\,\frac49(ab)^6\biggr)\!\biggr\},
%%% проверено в Mathematica %%%%%%%%%%%%%%%%%%%%%%%%%%%%%%%%%%%%%%%%%%%%%%
\label{eq:GaussAiBi_Integ1}
\end{align}
which is valid for all $a,b\in\Cset$.

% which is valid for all $a,b\in\Cset$ due to a standard theorem in complex 
% analysis \cite[Sect.\,2.6]{Titchmarsh_TheorFunc_ENG}.

% (two analytic functions which are equal in a region having a limit-point 
% inside, are equal identically in a region of joint analyticity).

In particular, if $b=-a$, then the integrand is the even function and
\begin{align}
\int_\Rset &\;\rme^{-t^2}\Ai(2at+a^4)\Bi(-2at+a^4)\,\rmd t
\nonumber\\
&{}=\frac{1}{12^{1/3}}\biggl\{\frac{1}{\Gamma\bigl(\tfrac56\bigr)}\cdot
{}_0F_1\biggl(\frac56\,\Bigl|\,\frac49a^{12}\biggr)
+\frac{2a^2}{12^{1/3}\Gamma\bigl(\tfrac76\bigr)}\cdot
{}_0F_1\biggl(\frac76\,\Bigl|\,\frac49a^{12}\biggr)\!\biggr\}.
%%% проверено в Mathematica %%%%%%%%%%%%%%%%%%%%%%%%%%%%%%%%%%%%%%%%%%%%%%
\label{eq:GaussAiBi_Integ2}
\end{align}
In the same way as in (\ref{eq:ff_fg_gg_Series}), we can use Airy atoms and 
separate the relation (\ref{eq:GaussAiBi_Integ2}) into three parts:
\begin{gather}
\int_\Rset \rme^{-t^2}f(2at+a^4)f(-2at+a^4)\,\rmd t
=\sqrt{\pi}\cdot{}_0F_1\biggl(\frac56\,\Bigl|\,\frac49a^{12}\biggr),
\nonumber\\
\int_\Rset \rme^{-t^2}g(2at+a^4)g(-2at+a^4)\,\rmd t
=-2\sqrt{\pi}\,a^2\cdot{}_0F_1\biggl(\frac76\,\Bigl|\,\frac49a^{12}\biggr),
\nonumber\\
\int_\Rset \rme^{-t^2}f(2at+a^4)g(-2at+a^4)\,\rmd t
=\int_\Rset \rme^{-t^2}g(2at+a^4)f(-2at+a^4)\,\rmd t=0.
\label{eq:Integ_AiryAtoms}
\end{gather}
Combining them properly, it is easy to find similar integrals for other 
products of Airy functions:
\begin{align}
\int_\Rset \rme^{-t^2}\Bi(2at+a^4)\Bi(-2at+a^4)\,\rmd t
=3\!\int_\Rset \rme^{-t^2}\Ai(2at+a^4)\Ai(-2at+a^4)\,\rmd t
\nonumber\\
{}=\frac{\sqrt 3}{12^{1/3}}\biggl\{\frac{1}{\Gamma\bigl(\tfrac56\bigr)}\cdot
{}_0F_1\biggl(\frac56\,\Bigl|\,\frac49a^{12}\biggr)
-\frac{2a^2}{12^{1/3}\Gamma\bigl(\tfrac76\bigr)}\cdot
{}_0F_1\biggl(\frac76\,\Bigl|\,\frac49a^{12}\biggr)\!\biggr\}.
%%% проверено в Mathematica %%%%%%%%%%%%%%%%%%%%%%%%%%%%%%%%%%%%%%%%%%%%%%
\label{eq:GaussAiBi_Integ3}
\end{align}

For other two cases, $E(z^2,0,z)$ and $E(-z^2/3,0,z)$, the relations are as 
follows:
\begin{align}
\int_\Rset \rme^{-t^2}&\bigl\{\Ai(t_1)\Bi(t_2)+\Ai(t_2)\Bi(t_1)\bigr\}\,\rmd t
=\frac{2}{12^{1/3}}\exp\biggl(-\frac13(a-b)^4(a^2+4ab+b^2)\biggr)
\nonumber\\
&{}\times\biggl\{\frac{1}{\Gamma\bigl(\tfrac56\bigr)}\cdot
{}_2F_2\biggl(\!\!\begin{array}{c} \frac{1}{12},\frac{7}{12} \smallskip\\ 
\frac13,\frac23\end{array}\!\Bigl|\,\frac{64}{3}(ab)^3\biggr)
-\frac{2ab}{12^{1/3}\Gamma\bigl(\tfrac76\bigr)}\cdot
{}_2F_2\biggl(\!\!\begin{array}{c} \frac{5}{12},\frac{11}{12} \smallskip\\ 
\frac23,\frac43\end{array}\!\Bigl|\,\frac{64}{3}(ab)^3\biggr)\!\biggr\},
\label{eq:GaussAiBi_Integ4}
\end{align}
where $t_1=2at-a(a-b)^2(a+2b)$ and $t_2=2bt-b(a-b)^2(2a+b)$,
\begin{align}
\int_\Rset \rme^{-t^2}&\bigl\{\Ai(t_3)\Bi(t_4)+\Ai(t_4)\Bi(t_3)\bigr\}\,\rmd t
=\frac{2}{12^{1/3}}
\exp\biggl(-\frac13\bigl([a^3+b^3]^2+3a^2b^2[a+b]^2\bigr)\biggr)
\nonumber\\
&{}\times\biggl\{\frac{1}{\Gamma\bigl(\tfrac56\bigr)}\cdot
{}_2F_2\biggl(\!\!\begin{array}{c} \frac16,\frac12 \smallskip\\ 
\frac13,\frac23\end{array}\!\Bigl|\,\frac{16}{3}(ab)^3\biggr)
-\frac{2ab}{12^{1/3}\Gamma\bigl(\tfrac76\bigr)}\cdot
{}_2F_2\biggl(\!\!\begin{array}{c} \frac12,\frac56 \smallskip\\ 
\frac23,\frac43\end{array}\!\Bigl|\,\frac{16}{3}(ab)^3\biggr)\!\biggr\},
%%% проверено в Mathematica %%%%%%%%%%%%%%%%%%%%%%%%%%%%%%%%%%%%%%%%%%%%%%
\label{eq:GaussAiBi_Integ5}
\end{align}
where $t_3=2at-a(a+b)(a^2-ab+2b^2)$ and $t_4=2bt-b(a+b)(2a^2-ab+b^2)$.
The simplest corollaries of both formulae (when $b=a$ and $b=-a$, 
respectively) may be proven directly, without using the theory above.

\section{Concluding remarks}
\label{sec:Extro}
%%%%%%%%%%%%%%%%%%%%%%%%%%%%%%%%%%%%%%%%%%%%%%%%%%%%%%%%%%%%%%%%%%%%%%%%%%%%

In the paper, we obtain the absolutely convergent integral representation 
for the elliptic umbilic and derive some of its corollaries. The formulae 
presented above may be useful for establishing possibly new hypergeometric 
identities and relations containing Airy polynomials. We discuss them briefly.

Since $\omega^{n-m+1}=\omega^{1-m-2n}$, we can transform 
(\ref{eq:Umbilic_XRc}) and (\ref{eq:Umbilic_XLc}) collecting terms 
with $m+2n=k$. It leads to the elliptic umbilic expansions which are valid 
for all $c$. They contain ${}_2F_2$ hypergeometric functions but have a more 
complicated structure than (\ref{eq:Umbilic_CausticXR}) and 
(\ref{eq:Umbilic_CausticXL}).

The hypergeometric functions appear in other expansions of $E(x,y,z)$, 
presented in Section~\ref{sec:IntegAiry}. For example, equating the 
coefficients of $z^n$ in (\ref{eq:Umbilic_HyperG_00z}) and 
(\ref{eq:Umbilic_Series-z/3}), where we put $x=y=0$, gives a closed-form 
expression for ${}_2F_1\bigl(a,a+\tfrac13;\,a+\tfrac23\,\bigl|\,
\tfrac12\bigr)$, which is known. (More exactly, we first obtain the value of 
${}_2F_1$ function for discrete $n$ and then replace $n$ by a continuous 
parameter, $-n/3\to a$. The last step can be justified by Carlson's theorem 
\cite[\S\,5.8.1]{Titchmarsh_TheorFunc_ENG}; see also \cite[Appendix]{Sigma}.) 

Applying this approach to (\ref{eq:Umbilic_HyperG_00z}) and 
(\ref{eq:Umbilic_Series+z/3}), we find that
\begin{equation}
{}_3F_2\biggl(\!\begin{array}{c} a,a+\tfrac12,\tfrac{4a+1}{6}\smallskip\\
a+\tfrac14,a+\tfrac34 \end{array}\!\Bigl|\,\frac34\biggr)
=\frac{\Gamma\bigl(\tfrac{5-4a}{6}\bigr)
\Gamma\bigl(\tfrac{1+2a}{6}\bigr)\Gamma\bigl(1+\tfrac{a}{3}\bigr)}
{2^{2-2a}\sqrt{\pi}\,\Gamma\bigl(\tfrac12-2a\bigr)\Gamma(1+2a)
\cos\tfrac{2\pi a}{3}\cos\tfrac{\pi(1+2a)}{3}}.
\label{eq:HyperG_3F2_x=3/4_n2}
%%% проверено в Mathematica %%%%%%%%%%%%%%%%%%%%%%%%%%%%%%%%%%%%%%%%%%%%%%
\end{equation}
Next, equating the coefficients of $x^\ell y^m z^n$ in 
(\ref{eq:Umbilic_Series-z/3}) and (\ref{eq:Umbilic_Series+z/3}) leads to
the value of ${}_3F_2\bigl(*\,\bigl|\,\tfrac34\bigr)$ depending on three 
parameters, which is too cumbersome. Besides, replacing discrete parameters 
$\ell,m,n$ by continuous ones, we could not check the conditions of Carlson's 
theorem.

Finally, Eq.\,(\ref{eq:Umbilic4}) provides a way to connect the polynomials 
$P_n(x)$, $Q_n(x)$ in (\ref{eq:AiryDeriv}) with Airy polynomials 
\cite{Torre2012}. The last ones are a particular case of the Kamp\'e de 
F\'eriet polynomials \cite{Dattoli} and may be defined by the generating 
function
\begin{equation}
\sum_{n\ge 0} \mathrm{Pi}_n(x)\frac{t^n}{n!}
=\exp\biggl(xt-\frac{t^3}{3}\biggr).
\label{eq:GenFunc_Pi}
\end{equation}
Substituting $a=\rmi c\sqrt\gamma$ and $b=0$ into (\ref{eq:Umbilic4}), 
after some algebra we have
\begin{equation}
\sum_{n\ge 0} \,\Ai^{(2n)}(x)\,\frac{c^n}{n!}
=\exp\biggl(cx+\frac{2c^3}{3}\biggr)\Ai(x+c^2).
\label{eq:Series_AiDeriv}
\end{equation}
Then
\begin{equation}
\Ai^{(2n)}(x)=n!\cdot\lshad c^n\rshad\biggl\{
\sum_{\ell\ge 0} \mathrm{Pi}_\ell\biggl(-\frac{x}{2^{1/3}}\biggr)
\frac{(-2^{1/3}c)^\ell}{\ell!}
\sum_{k\ge 0} \,\Ai^{(k)}(x)\,\frac{c^{2k}}{k!}\biggr\},
\label{eq:Series_AiDeriv1}
\end{equation}
where we use the notation proposed in \cite{GKP_ENG}. Namely, if 
$A(z)$ is any power series $\sum_k a_kz^k$, then $\lshad z^k\rshad A(z)$ 
denotes the coefficient of $z^k$ in $A(z)$. In our view, this notation is 
more convenient to manipulate power series than usual analytic description, 
$\lshad z^k\rshad A(z)=A^{(k)}(0)/k!$.

Separating the components of (\ref{eq:Series_AiDeriv1}), one gets the 
relations 
\begin{equation}
\begin{aligned}
P_{2n}(x)&=\sum_{0\le 2k\le n} \frac{n!\,(-2^{1/3})^{n-2k}}{(n-2k)!\,k!}
\,\mathrm{Pi}_{n-2k}\biggl(-\frac{x}{2^{1/3}}\biggr)P_k(x),\\
Q_{2n}(x)&=\sum_{0\le 2k\le n} \frac{n!\,(-2^{1/3})^{n-2k}}{(n-2k)!\,k!}
\,\mathrm{Pi}_{n-2k}\biggl(-\frac{x}{2^{1/3}}\biggr)Q_k(x).
\end{aligned}
\label{eq:PQ_Pi}
\end{equation}

\appendix{}
\section{The power series for $\Re\{\Ai(x+\rmi y)\Bi(x-\rmi y)\}$}
\label{sec:AppA}
\setcounter{equation}{0}
\renewcommand{\theequation}{A{\arabic{equation}}}
%%%%%%%%%%%%%%%%%%%%%%%%%%%%%%%%%%%%%%%%%%%%%%%%%%%%%%%%%%%%%%%%%%%%%%%%%%%%

Let us find the expansion of $\Re\{\Ai(x+\rmi y)\Bi(x-\rmi y)\}$ in power 
series. We deduce it from the expansion for $\Ai(x+y)\Ai(x-y)$. There are 
at least two ways to prove that
\begin{equation}
\Ai(x+y)\Ai(x-y)=\frac{2}{\sqrt\pi}
\sum_{m,n\ge 0} \;\frac{(-1)^{m+n}x^m y^{2n}}
{12^{(5+2n-2m)/6}\Gamma\bigl(\tfrac{5+2n-2m}{6}\bigr)\,m!\,n!}.
%%% проверено в Mathematica %%%%%%%%%%%%%%%%%%%%%%%%%%%%%%%%%%%%%%%%%%%%%%
\label{eq:AiAi_Series}
\end{equation}
Here, we use the approach based on the Mellin transform,
\begin{equation}
\Mgoth[f(x)](\alpha)=\int_0^\infty x^{\alpha-1}f(x)\,\rmd x,
\qquad \alpha>0,
\label{eq:Mellin}
\end{equation}
applying to the shifted Airy functions.

We begin with the well-known formula \cite[Eq.\,(2.21)]{Reid}
\begin{equation}
\Ai^2(x)+\Bi^2(x)=\frac{1}{\pi^{3/2}}
\int_0^\infty \exp\biggl(xt-\frac{t^3}{12}\biggr)\frac{\rmd t}{\sqrt t},
\label{eq:Reid_Ai2+Bi2}
\end{equation}
replacing $x$ by $\omega x$ and using (\ref{eq:AlgRels}):
\begin{gather}
\Ai^2(\omega x)+\Bi^2(\omega x)
=2\rme^{-\pi\rmi/3}\bigl\{\Ai^2(x)+\rmi\cdot\Ai(x)\Bi(x)\bigr\},
\label{eq:Ai2+Bi2=Ai2+AiBi}\\
\Ai^2(x)+\rmi\cdot\Ai(x)\Bi(x)=\frac{\rme^{\pi\rmi/3}}{2\pi^{3/2}}
\int_0^\infty \exp\biggl(\omega xt-\frac{t^3}{12}\biggr)\frac{\rmd t}{\sqrt t}.
\label{eq:Integ1_Ai^2+iAiBi}
\end{gather}
Then
\begin{align}
\int_0^\infty &{}x^{\alpha-1}\bigl\{\Ai^2+\rmi\cdot\Ai\Bi\bigr\}(x+c)\,\rmd x
\nonumber\\
&{}=\frac{\rme^{\pi\rmi/3}}{2\pi^{3/2}}
\int_0^\infty \exp\biggl(\omega ct-\frac{t^3}{12}\biggr)\frac{\rmd t}{\sqrt t}
\int_0^\infty \exp\bigl(-\rme^{-\pi\rmi/3}tx\bigr) x^{\alpha-1}\,\rmd x
\nonumber\\
&{}=\frac{\rme^{\pi\rmi/3}}{2\pi^{3/2}}
\int_0^\infty \exp\biggl(\omega ct-\frac{t^3}{12}\biggr)\cdot
\frac{\Gamma(\alpha)}{(\rme^{-\pi\rmi/3}t)^\alpha}\cdot\frac{\rmd t}{\sqrt t}
\nonumber\\
&{}=\frac{\Gamma(\alpha)}{2\pi^{3/2}}\,\rme^{\pi\rmi(\alpha+1)/3}
\sum_{n\ge 0} \frac{(\omega c)^n}{n!}
\int_0^\infty \exp\biggl(-\frac{t^3}{12}\biggr)t^{n-\alpha-1/2}\,\rmd t
\nonumber\\
&{}=\frac{\Gamma(\alpha)}{6\pi^{3/2}}\,\rme^{\pi\rmi(\alpha+1)/3}
\sum_{n\ge 0} \Gamma\biggl(\frac{2n+1-2\alpha}{6}\biggr)
12^{(2n+1-2\alpha)/6}\frac{(\omega c)^n}{n!}
\nonumber\\
&{}=\frac{2\Gamma(\alpha)}{\sqrt\pi}\sum_{n\ge 0} 
\,\frac{1-\rmi\cdot\cot\bigl(\pi\,\tfrac{5+2\alpha-2n}{6}\bigr)}
{12^{(5+2\alpha-2n)/6}\Gamma\bigl(\tfrac{5+2\alpha-2n}{6}\bigr)}
\cdot\frac{(-c)^n}{n!}.
\label{eq:Ai^2+AiBic_Mellin}
\end{align}
The last equality follows after using the reflection formula: 
$\Gamma(z)\Gamma(1-z)=\pi/\sin(\pi z)$.
To insure the convergence of the integral we restrict our study to the 
case of small $\alpha$’s, namely, $\alpha\in\bigl(0,\tfrac12\bigr)$. 
Of course, this restriction is due to the imaginary part only. As for the 
real part, the resulting expression 
\begin{equation}
\Mgoth\bigl[\Ai^2(x+c)\bigr](\alpha)
=\frac{2\Gamma(\alpha)}{\sqrt\pi}\sum_{n\ge 0} 
\,\frac{12^{(2n-2\alpha-5)/6}}{\Gamma\bigl(\tfrac{5+2\alpha-2n}{6}\bigr)}
\cdot\frac{(-c)^n}{n!}
\label{eq:Shifted_Ai^2_Mellin}
\end{equation}
is valid for all $\alpha>0$ and may be established by analytic continuation 
(see \cite{Watson_ENG,Nikiforov_ENG} for details).

Next, we use the formula \cite[Eq.\,(B18)]{Berry1979}
\begin{equation}
\Ai\biggl(\frac{a+b}{2^{2/3}}\biggr)\Ai\biggl(\frac{a-b}{2^{2/3}}\biggr)
=\frac{1}{2^{2/3}\pi}\int_\Rset \rme^{\rmi bt}\Ai(a+t^2)\,\rmd t
\label{eq:B18}
\end{equation}
to get the relation between the Mellin transforms of the shifted $\Ai$ and 
$\Ai^2$ functions:
\begin{align}
\int_0^\infty &{}x^{\alpha-1}\Ai\biggl(\frac{x+c+b}{2^{2/3}}\biggr)
\Ai\biggl(\frac{x+c-b}{2^{2/3}}\biggr)\,\rmd x
\nonumber\\
&{}=\frac{1}{2^{2/3}\pi}\int_0^\infty x^{\alpha-1}\,\rmd x
\int_\Rset \rme^{\rmi bt}\Ai(x+c+t^2)\,\rmd t
\nonumber\\
&{}=\frac{1}{2^{2/3}\pi}
\int_0^\infty \bigl(\rme^{\rmi bt}+\rme^{-\rmi bt}\bigr)\,\rmd t
\int_0^\infty x^{\alpha-1}\Ai(x+t^2+c)\,\rmd x
\nonumber\\
&{}=\frac{4}{2^{2/3}\pi}\int_0^\infty \!\!\int_0^\infty 
\cos(bt)\,s^{2\alpha-1}\Ai(s^2+t^2+c)\,\rmd s\,\rmd t
\nonumber\\
&{}=\frac{4}{2^{2/3}\pi}\int_0^\infty r^{2\alpha}\Ai(r^2+c)\,\rmd r
\int_0^{\pi/2} \cos(br\sin\phi)\cos^{2\alpha-1}\phi\,\rmd\phi
\nonumber\\
&{}=\frac{4}{2^{2/3}\pi}\int_0^\infty r^{2\alpha}\Ai(r^2+c)
\sum_{n\ge 0} \frac{(-b^2r^2)^n}{(2n)!}\cdot\frac12\,
\Beta\biggl(n+\frac12,\alpha\biggr)\,\rmd r
\nonumber\\
&{}=\frac{2\Gamma(\alpha)}{2^{2/3}\pi}\sum_{n\ge 0} 
\;\frac{\Gamma\bigl(n+\tfrac12\bigr)(-b^2)^n}
{\Gamma\bigl(n+\alpha+\tfrac12\bigr)(2n)!}\,
\int_0^\infty r^{2n+2\alpha}\Ai(r^2+c)\,\rmd r
\nonumber\\
&{}=\frac{2\Gamma(\alpha)}{2^{2/3}\sqrt\pi}\sum_{n\ge 0} 
\;\frac{\bigl(-\tfrac14 b^2\bigr)^n}{\Gamma\bigl(n+\alpha+\tfrac12\bigr)\,n!}\,
\int_0^\infty r^{2n+2\alpha}\Ai(r^2+c)\,\rmd r
\nonumber\\
&{}=\frac{2\Gamma(\alpha)}
{2^{2/3}\sqrt{\pi}\,\Gamma\bigl(\alpha+\tfrac12\bigr)}\,\int_0^\infty 
{}_0F_1\biggl(\alpha+\frac12\,\Bigl|\,-\frac{b^2r^2}{4}\biggr)\, 
r^{2\alpha}\Ai(r^2+c)\,\rmd r
\nonumber\\
&{}=\frac{\Gamma(\alpha)}
{2^{2/3}\sqrt{\pi}\,\Gamma\bigl(\alpha+\tfrac12\bigr)}\,\int_0^\infty 
{}_0F_1\biggl(\alpha+\frac12\,\Bigl|\,-\frac{b^2x}{4}\biggr)\, 
x^{\alpha-1/2}\Ai(x+c)\,\rmd x.
\label{eq:AiAi_Mellin}
\end{align}
Neglecting quite interesting integral relations between Airy and Bessel 
functions (see also \cite{Varlamov2010}), we substitute $b=0$ and obtain
\begin{equation}
2^{2\alpha/3}\Mgoth\biggl[\Ai^2\biggl(x+\frac{c}{2^{2/3}}\biggr)\biggr](\alpha)
=\frac{\Gamma(\alpha)}{2^{2/3}\sqrt{\pi}\,\Gamma\bigl(\alpha+\tfrac12\bigr)}\,
\Mgoth\bigl[\Ai(x+c)\bigr]\biggl(\alpha+\frac12\biggr).
\label{eq:Mellin_Ai2=Ai}
\end{equation}
As a corollary,
\begin{align}
\Mgoth\bigl[\Ai(x+c)\bigr](\alpha)
&=\frac{2^{(2\alpha+1)/3}\sqrt{\pi}\,\Gamma(\alpha)}
{\Gamma\bigl(\alpha-\tfrac12\bigr)}\,
\Mgoth\biggl[\Ai^2\biggl(x+\frac{c}{2^{2/3}}\biggr)\biggr]
\biggl(\alpha-\frac12\biggr)
\nonumber\\
&\overset{\makebox[0pt]{\scriptsize (\ref{eq:Shifted_Ai^2_Mellin})}}{=}
\,\Gamma(\alpha)\sum_{n\ge 0} 
\,\frac{3^{(n-\alpha-2)/3}}{\Gamma\bigl(\tfrac{\alpha+2-n}{3}\bigr)}
\cdot\frac{(-c)^n}{n!}.
\label{eq:Shifted_Ai_Mellin}
\end{align}

Now, one can find the power series for the product of two Airy functions:
\begin{align}
\Ai\biggl(\frac{a+b}{2^{2/3}}\biggr)\Ai\biggl(\frac{a-b}{2^{2/3}}\biggr)
\,&\overset{\makebox[0pt]{\scriptsize (\ref{eq:B18})}}{=}
\,\frac{1}{2^{2/3}\pi}\int_\Rset \rme^{\rmi bx}\Ai(a+x^2)\,\rmd x
\nonumber\\
&{}=\frac{1}{2^{2/3}\pi}\sum_{n\ge 0} \frac{(\rmi b)^n}{n!}
\int_\Rset x^n\Ai(a+x^2)\,\rmd x={}
\nonumber\\
&{}=\frac{1}{2^{2/3}\pi}\sum_{n\ge 0} \frac{(-b^2)^n}{(2n)!}
\int_0^\infty t^{n-1/2}\Ai(a+t)\,\rmd t
\nonumber\\
\,&\overset{\makebox[0pt]{\scriptsize (\ref{eq:Shifted_Ai_Mellin})}}{=}
\frac{1}{\sqrt\pi}\sum_{m\ge 0} \sum_{n\ge 0} 
\frac{3^{(m-5/2-n)/3}(-a)^m(-b^2)^n}
{2^{2n+2/3}\Gamma\bigl(\tfrac{5+2n-2m}{6}\bigr)\,m!\,n!}
\nonumber
\end{align}
which is analogous to the formula (\ref{eq:AiAi_Series}).

Because of the symmetry relations (\ref{eq:AlgRels}), the following formulae 
hold:
\begin{equation}
\begin{aligned}
4\Re\bigl\{\omega\Ai(\omega^*[x+y])\Ai(\omega^*[x-y])\bigr\}
&=\Ai(x+y)\Ai(x-y)-\Bi(x+y)\Bi(x-y),\\
4\Im\bigl\{\omega\Ai(\omega^*[x+y])\Ai(\omega^*[x-y])\bigr\}
&=\Ai(x+y)\Bi(x-y)+\Ai(x-y)\Bi(x+y).
\end{aligned}
\label{eq:AlgRels2}
\end{equation}
Then it is easy to get the power series for other products of Airy functions:
\begin{align}
\Ai(x+y)\Ai(x-y)&=\sum_{m,n\ge 0} c_{m,n}\quad\Rightarrow
\nonumber\\
\Bi(x+y)\Bi(x-y)&=\sum_{m,n\ge 0} \bigl\{1-4\Re\omega^{n-m+1}\bigr\}c_{m,n},
\nonumber\\
\Ai(x+y)\Bi(x-y)+\Ai(x-y)\Bi(x+y)
&=4\sum_{m,n\ge 0} \Im\{\omega^{n-m+1}\}\,c_{m,n},
\label{eq:BiBi_Series}
\end{align}
where we rewrite (\ref{eq:AiAi_Series}) in short-hand. The last expression 
helps to find the desired result:
\begin{equation}
\Re\{\Ai(x+\rmi y)\Bi(x-\rmi y)\}=\frac{4}{\sqrt\pi}
\sum_{m,n\ge 0} \;\frac{\Im\{\omega^{n-m+1}\}(-x)^m y^{2n}}
{12^{(5+2n-2m)/6}\Gamma\bigl(\tfrac{5+2n-2m}{6}\bigr)\,m!\,n!}.
\label{eq:ReAiBi}
\end{equation}

An alternative way is to use \textit{Airy atoms} $f(x)$ and $g(x)$ which are 
connected with Airy functions by the relations \cite{Abramowitz_ENG}:
\begin{equation}
\Ai(x)=c_1f(x)-c_2g(x),\qquad \frac{\Bi(x)}{\sqrt 3}=c_1f(x)+c_2g(x),
\label{eq:AiryAtoms1}
\end{equation}
where $c_1=3^{-2/3}\!/\Gamma\bigl(\tfrac23\bigr)$ and
$c_2=3^{-1/3}\!/\Gamma\bigl(\tfrac13\bigr)$. Since the coefficients of power 
series for both Airy atoms are \textit{rational},
\begin{equation}
\begin{aligned}
f(x)&=\sum_{k\ge 0} \left(\frac13\right)_{\!k}\frac{3^k x^{3k}}{(3k)!}
={}_0F_1\left(\frac23\,\Bigl|\,\frac{x^3}{9}\right),\\
g(x)&=\sum_{k\ge 0} \left(\frac23\right)_{\!k}\frac{3^k x^{3k+1}}{(3k+1)!}
=x\cdot{}_0F_1\left(\frac43\,\Bigl|\,\frac{x^3}{9}\right),
\end{aligned}
\label{eq:AiryAtoms}
\end{equation}
we can separate (\ref{eq:AiAi_Series}) into three parts:
\begin{gather}
\begin{aligned}
&\Ai(x+y)\Ai(x-y)=\{c_1 f(x+y)-c_2 g(x+y)\}\{c_1 f(x-y)-c_2 g(x-y)\}\\
&{}=\frac{f(x+y)f(x-y)}{12^{1/3}\sqrt{3\pi}\,\Gamma\bigl(\tfrac56\bigr)} 
-\frac{f(x-y)g(x+y)+f(x+y)g(x-y)}{2\pi\sqrt{3}}
+\frac{g(x+y)g(x-y)}{12^{1/6}\sqrt{\pi}\,\Gamma\bigl(\tfrac16\bigr)},
\end{aligned}
\nonumber\\
\begin{gathered}
\frac{f(x+y)f(x-y)}{12^{1/3}\sqrt{3\pi}\,\Gamma\bigl(\tfrac56\bigr)}
=\sum_{\substack{m,n\ge 0\\ m-n\equiv 0(\mathrm{mod}\,3)}} c_{m,n},\\
-\frac{f(x-y)g(x+y)+f(x+y)g(x-y)}{2\pi\sqrt{3}}
=\sum_{\substack{m,n\ge 0\\ m-n\equiv 1(\mathrm{mod}\,3)}} c_{m,n},\\
\frac{g(x+y)g(x-y)}{12^{1/6}\sqrt{\pi}\,\Gamma\bigl(\tfrac16\bigr)}
=\sum_{\substack{m,n\ge 0\\ m-n\equiv 2(\mathrm{mod}\,3)}} c_{m,n}.
\end{gathered}
\label{eq:ff_fg_gg_Series}
\end{gather}
Applying these expansions to $\Re\{\Ai(x+\rmi y)\Bi(x-\rmi y)\}$, we get a 
double series
\begin{align}
&\Re\{\Ai(x+\rmi y)\Bi(x-\rmi y)\}
=\sqrt{3}\,c_1^2\bigl|f(x+\rmi y)\bigr|^2
-\sqrt{3}\,c_2^2\bigl|g(x+\rmi y)\bigr|^2
\nonumber\\
&=\frac{2\sqrt{3}}{\sqrt\pi}\,
\biggl\{\sum_{\substack{m,n\ge 0\\ m-n\equiv 0(\mathrm{mod}\,3)}}
-\sum_{\substack{m,n\ge 0\\ m-n\equiv 2(\mathrm{mod}\,3)}}\biggr\}
\,\frac{(-x)^m y^{2n}}
{12^{(5+2n-2m)/6}\Gamma\bigl(\tfrac{5+2n-2m}{6}\bigr)\,m!\,n!},
\label{eq:ReAiBi_2}
\end{align}
which evidently coincides with (\ref{eq:ReAiBi}).

As a final remark, we provide the counterpart of (\ref{eq:ReAiBi}) for 
completeness,
\begin{align}
\Im\{\Ai(x+\rmi y)\Bi(x-\rmi y)\}&=-\frac{1}{\pi}
\sum_{m,n\ge 0} 
\;\frac{(-x)^m y^{6n+2m+1}}{12^n\bigl(\tfrac32\bigr)_{3n+m}\,m!\,n!}
\nonumber\\
&=-\frac{1}{\pi}\sum_{k\ge 0} \;\frac{(-4)^k y^{2k+1}}{(2k+1)!}\cdot
\frac{\mathrm{Pi}_k(2^{2/3}x)}{(2^{2/3})^k},
\label{eq:ImAiBi}
\end{align}
and mention that another way to prove (\ref{eq:AiAi_Series}) 
is based on a slightly modified version of Moyer's formula \cite{Moyer} 
(see also \cite{Berry2013}),
\begin{equation}
\Ai(x+y)\bigl\{\Ai(x-y)+\rmi\cdot\Bi(x-y)\bigr\}
=\frac{1}{2\pi^{3/2}}\int_0^\infty 
\exp\biggl(\rmi\biggl[\frac{t^3}{12}+xt-\frac{y^2}{t}+\frac{\pi}{4}\biggr]
\biggr)\frac{\rmd t}{\sqrt t},
\label{eq:Moyer}
\end{equation}
namely,
\begin{equation}
\bigl|\Ai(x+\rmi y)+\rmi\cdot\Bi(x+\rmi y)\bigr|^2
=\frac{1}{\pi^{3/2}}\int_0^\infty 
\exp\biggl(-\frac{t^3}{12}+xt-\frac{y^2}{t}\biggr)\frac{\rmd t}{\sqrt t},
\label{eq:Moyer2}
\end{equation}
where $x\in\Rset$ and $y\ge 0$.

%\section*{References}
\label{sec:Refs}
%%%%%%%%%%%%%%%%%%%%%%%%%%%%%%%%%%%%%%%%%%%%%%%%%%%%%%%%%%%%%%%%%%%%%%%%%%%%

%\bibliographystyle{unsrt}
%\bibliography{all}

\end{document}